\documentclass[12pt]{amsart}
\usepackage{amsmath, amsfonts, amssymb, amsthm,hyperref,mathtools,array}
\hypersetup{hypertex=true,
	colorlinks=true,
	linkcolor=blue,
	anchorcolor=blue,
	citecolor=blue}
\usepackage{bm}
\allowdisplaybreaks[4]
\def\({\bg(}
\def\){\bg)}

\def\ord{{\rm ord}}

\def\Gal{{\rm Gal}}
\def\Tr{{\rm Tr}}

\def\0{{\bm 0}}

\def\diag{{\rm diag}}

\def\pmod #1{\ ({\rm{mod}}\ #1)}
\def\mod #1{\ {\rm mod}\ #1}
\def\Ack{\medskip\noindent {\bf Acknowledgments}}

\theoremstyle{plain}
\newtheorem{theorem}{Theorem}[section]
\newtheorem{lemma}{Lemma}

\theoremstyle{definition}

\theoremstyle{remark}
\newtheorem{remark}{Remark}

\makeatletter
\@namedef{subjclassname@2020}{%
	\textup{2020} Mathematics Subject Classification}
\makeatother
\vspace{4mm}

\begin{document}
	
	\title[Cyclotomic matrices involving Gauss sums over cyclic groups]
	{The Gauss periods and cyclotomic matrices involving Gauss sums over cyclic groups}
	\author[X.-H. Ji, L.-Y. Wang and H.-L. Wu]{Xiao-Han Ji, Li-Yuan Wang* and Hai-Liang Wu}
	
	\address {(Xiao-Han Ji) School of Science, Nanjing University of Posts and Telecommunications, Nanjing 210023, People's Republic of China}
	\email{\tt jxhmath@163.com}
	
   \address {(Li-Yuan Wang) School of Physical and Mathematical Sciences, Nanjing Tech University, Nanjing 211816, People's Republic of China}
   \email{\tt wly@smail.nju.edu.cn}
   
   \address {(Hai-Liang Wu) School of Science, Nanjing University of Posts and Telecommunications, Nanjing 210023, People's Republic of China}
   \email{\tt whl.math@smail.nju.edu.cn}

	\keywords{Gauss sums, cyclotomic matrices, finite fields.
		\newline \indent This research was supported by the Natural Science Foundation of China (Grant Nos. 12671009 and 12201291) and the Natural Science Foundation of the Higher Education Institutions of Jiangsu Province (Grant No. 25KJB110010).
		\newline \indent *Corresponding author.}
	
	\begin{abstract}
		In this paper, by using the arithmetic properties of the Gauss periods and character sums over cyclic groups, we study the cyclotomic matrix 
		$$A_k(\chi)=\left[G_N(\chi^{ki+ki})\right]_{0\le i,j\le \varphi(N)/k-1},$$
		where $N=p^m$ is an odd prime power, $\varphi(\cdot)$ is the Euler totient function, $k$ is a divisor of $\varphi(N)$, $\chi$ is a generator of character group $\widehat{(\mathbb{Z}/N\mathbb{Z})^{\times}}$, and 
		$$G_N(\chi^{ki+kj})=\sum_{x\in\mathbb{Z}/N\mathbb{Z}}\chi^{ki+kj}(x)e^{2\pi ix/N}$$
		is the Gauss sum over $\mathbb{Z}/N\mathbb{Z}$.
	\end{abstract}
	\maketitle
	
	\section{Introduction}
	\setcounter{lemma}{0}
	\setcounter{theorem}{0}
	\setcounter{equation}{0}
	\setcounter{conjecture}{0}
	\setcounter{remark}{0}
	\setcounter{corollary}{0}
	
	\subsection{Notation} For any integer $N>1$, we use $\mathbb{Z}/N\mathbb{Z}$ to denote the ring of residue classes modulo $N$, and let 
	$$\left(\mathbb{Z}/N\mathbb{Z}\right)^{\times}=\left\{x\mod{N\mathbb{Z}}: \gcd(x,N)=1\right\}$$
	be the multiplicative group of all invertible elements in $\mathbb{Z}/N\mathbb{Z}$. A group homomorphism 
	$$\psi: \left(\mathbb{Z}/N\mathbb{Z}\right)^{\times}\rightarrow\mathbb{C}^{\times}$$
	is called a Dirichlet character modulo $N$. The group of all Dirichlet characters modulo $N$ is known as the dual group of $(\mathbb{Z}/N\mathbb{Z})^{\times}$ and denoted by $\widehat{(\mathbb{Z}/N\mathbb{Z})^{\times}}$. In addition, given any $\psi\in\widehat{(\mathbb{Z}/N\mathbb{Z})^{\times}}$, we set $\psi(x\mod{N\mathbb{Z}})=0$ for any non-invertible element $x\mod{N\mathbb{Z}}\in\mathbb{Z}/N\mathbb{Z}$. 
	
	Let  $\zeta_N=e^{2\pi {\bm i}/N}$, where ${\bm i}$ is a primitive $4$th root of unity with argument $\pi/2$. Then, for any $\psi\in\widehat{(\mathbb{Z}/N\mathbb{Z})^{\times}}$, the Gauss sum  $G_N(\psi)$ over $\mathbb{Z}/N\mathbb{Z}$ is defined by 
	$$G_N(\psi)=\sum_{x\in\mathbb{Z}/N\mathbb{Z}}\psi(x)\zeta_N^x=\sum_{x\in(\mathbb{Z}/N\mathbb{Z})^{\times}}\psi(x)\zeta_N^x.$$
	
	Also, for any square matrix $M$ over a field, the symbol $M(i,j)$ denotes the $(i,j)$-entry of $M$, and $\det M$ indicates the determinant of $M$. 
	
	\subsection{Background and motivation} Let $p$ be an odd prime. In 1811, the exact value of quadratic Gauss sum over $\mathbb{Z}/p\mathbb{Z}=\mathbb{F}_p$ was successfully determined by Gauss, that is, 
	\begin{equation}\label{Eq. Gauss sums over Fp}
		\sum_{x\in\mathbb{F}_p}\left(\frac{x}{p}\right)\zeta_p^x=\sqrt{(-1)^{(p-1)/2}p},
	\end{equation}
	where $(\frac{\cdot}{p})$ is the Legendre symbol. The above result can be generalized in two directions. One is to an arbitrary finite field with odd characteristic, and the other is to the ring $\mathbb{Z}/N\mathbb{Z}$. For example, in the first direction, using the Hasse-Davenport lifting formula (cf. \cite[Theorem 3.7.4]{Cohen}), one can generalize (\ref{Eq. Gauss sums over Fp}) to the finite field $\mathbb{F}_q$ with $q$ elements, where $q=p^m$. Specifically, 
	\begin{equation*}
		\sum_{x\in\mathbb{F}_q}\phi_q(x)\zeta_p^{\Tr(x)}=(-1)^{m-1}\cdot {\bm i}^{\frac{m(p-1)^2}{4}}\cdot \sqrt{q},
	\end{equation*}
	where $\phi_q$ is the unique quadratic multiplicative character of $\mathbb{F}_q$, and $\Tr(\cdot)$ is the trace map from $\mathbb{F}_q$ to $\mathbb{F}_p$. For the second direction, let $N\ge 3$ be an odd integer and $(\frac{\cdot}{N})$ be the Jacobi symbol. Then, as another generalization of (\ref{Eq. Gauss sums over Fp}), we have (cf. \cite[Theorem 3.3]{Iwaniec})
	$$\sum_{x\in\mathbb{Z}/N\mathbb{Z}}\left(\frac{x}{N}\right)\zeta_N^x=\mu(N)^2\cdot {\bm i}^{(N-1)^2/4}\cdot \sqrt{N},$$
	where $\mu(\cdot)$ is the M\"obius function. 
	
	Next we introduce some related works concerning cyclotomic matrices. For any non-trivial Dirichlet character $\psi$ modulo an odd prime $p$, Carlitz \cite[Theorem 5]{Carlitz} showed that 
	$$\det \left[\psi(i+j)\right]_{1\le i,j\le p-1}=\begin{cases}
		(-1)^{(p-1)/(2\ord(\psi))}G_p(\psi)^{p-1}/p  &  \mbox{if}\ \ord(\psi)\equiv 1\pmod 2,\\
		(-1)^{(p-1)/(\ord(\psi))}\delta(\psi)^{p-1}G_p(\psi)^{p-1}/p  &  \mbox{if}\ \ord(\psi)\equiv 0\pmod 2,
	\end{cases}$$
	where $\ord(\psi)$ denotes the order of the character $\psi$ and 
	$$\delta(\psi)=\begin{cases}
		1                &   \mbox{if}\ \psi(-1)=1,\\
		-{\bm i}    &   \mbox{if}\ \psi(-1)=1.
	\end{cases}$$
	This result reveals a close connection between the determinants of certain cyclotomic matrices and Gauss sums. Compared with Carlitz's result, in 2025, the author, Li, Wang and Yip \cite{WLWY} considered the cyclotomic matrices whose entries are Gauss sums, and proved that 
	$$\det\left[G_p(\chi_p^{i+j})\right]_{0\le i,j\le p-2}=(-1)^{(p-3)/2}(p-1)^{p-1},$$
	where $\chi_p$ is a generator of $\widehat{\mathbb{F}_p^{\times}}$. More generally, given any positive divisor $k$ of $p-1$ with $p-1=kn$, recently the author, Wang and Pan \cite[Theorem 1.2]{WWP} showed that  
	\begin{equation}\label{Eq. result by WWP}
		\det\left[G_p(\chi_p^{ki+kj})\right]_{0\le i,j\le n-1}=(-1)^{(n^2-n+2)/2}\cdot n^n\cdot y_p(k),
	\end{equation}
	where $y_p(k)$ is the constant term of the minimal polynomial of the algebraic integer 
	\begin{equation}\label{Eq. definition of theta k}
		\theta_p(k)=\sum_{\substack{x\in\mathbb{F}_p\\ x^k=1}}\zeta_p^x.
	\end{equation}
	
	Motivated by the above results, it is natural to investigate cyclotomic matrices with Gauss sums over $\mathbb{Z}/p^m\mathbb{Z}$ as entries. In the remaining part of this paper, we always let $N=p^m$, $n=\varphi(p^m)=p^{m-1}(p-1)$ and $k\ge 1$ be a divisor of $n$ with $n=kd$, where $\varphi(\cdot)$ is the Euler totient function. Since $(\mathbb{Z}/N\mathbb{Z})^{\times}$ is a cyclic group of order $n$ and $(\mathbb{Z}/N\mathbb{Z})^{\times}\cong \widehat{(\mathbb{Z}/N\mathbb{Z})^{\times}}$, we let $\chi$ be a generator of $\widehat{(\mathbb{Z}/N\mathbb{Z})^{\times}}$. With the above notations, we focus on the cyclotomic matrix 
	$$A_k(\chi)=\left[G_N(\chi^{ki+ki})\right]_{0\le i,j\le d-1}.$$
	
	Although the matrices $A_k(\chi)$ and  $[G_p(\chi_p^{ki+kj})]_{0\le i,j\le n-1}$ look similar, the methods for handling them are completely different. This is mainly because 
	$\mathbb{Z}/N\mathbb{Z}$ is not a field whenever $m\ge2$, which prevents us from using tools related to finite fields. 
	
	We will see below that $\det A_k(\chi)$ has a close relationship with the Gauss periods. For this reason, we briefly introduce the Gauss periods here. Readers may refer to \cite{BEW} for a thorough introduction on this topic. Let notations be as above and let $H$ be a subgroup of $(\mathbb{Z}/N\mathbb{Z})^{\times}$. For any $x\in (\mathbb{Z}/N\mathbb{Z})^{\times}$, the sum 
	$$\sum_{y\in xH}\zeta_N^y$$
	is called a Gauss period. Note that in the case $m=1$, letting $H=\{x\in\mathbb{F}_p: x^k=1\}\leq\mathbb{F}_p^{\times}$, the algebraic integer $\theta_p(k)$ defined by (\ref*{Eq. definition of theta k}) is exactly a Gauss period.

	\subsection{Main results} Let notations be as above and let 
	$$H_d=\left\{x^d: x\in(\mathbb{Z}/N\mathbb{Z})^{\times}\right\}\leq (\mathbb{Z}/N\mathbb{Z})^{\times},$$
	where $N=p^m$ be an odd prime power. We set the Gauss period 
	$$\eta_N(k):=\sum_{x\in H_d}\zeta_N^x,$$
	and let $P_k(T)$ be the minimal polynomial of $\eta_N(k)$ over $\mathbb{Q}$. Now we state our main result. 
	
	\begin{theorem}\label{Thm. A}
		Let $N=p^m$ be an odd prime power with $p$ prime and $m\in\mathbb{Z}^+$. Let $n=\varphi(N)$ and $k\ge 1$ be a divisor of $n$ with $n=kd$. Then, for any generator $\chi$ of $\widehat{(\mathbb{Z}/N\mathbb{Z})^{\times}}$ the following results hold.
		
		{\rm (i)} The matrix $A_k(\chi)$ is singular if and only if $k\equiv 0\pmod p$. 
		
		{\rm (ii)} Suppose $k\not\equiv 0\pmod p$. Then $\mathbb{Q}(\eta_N(k))$ is the unique intermediate field of the cyclic Galois extension  $\mathbb{Q}(\zeta_N)/\mathbb{Q}$ such that  $[\mathbb{Q}(\eta_N(k)):\mathbb{Q}]=d$. Moreover, 
		$$\det A_k(\chi)=(-1)^{d+\lfloor\frac{d-1}{2}\rfloor}\cdot d^d\cdot y_N(k)\in\mathbb{Z},$$
		where $y_N(k)\in\mathbb{Z}$ is the constant term of the minimal polynomial $P_k(T)$.
	\end{theorem}
	
	\begin{remark}
		Since the Gauss period $\eta_N(k)$ is independent of the generator $\chi$ of $\widehat{(\mathbb{Z}/N\mathbb{Z})^{\times}}$, it follows from Theorem \ref{Thm. A} that $\det A_k(\chi)$ is also independent of the choice of $\chi$. 
	\end{remark}
	
	\subsection{Outline of this paper} We will prove Theorem \ref{Thm. A} in Section 2.

	\section{Proof of Theorem \ref{Thm. A}}
	\setcounter{lemma}{0}
	\setcounter{theorem}{0}
	\setcounter{equation}{0}
	\setcounter{conjecture}{0}
	\setcounter{remark}{0}
	\setcounter{corollary}{0}
	
	For any positive integer $m$, we let $\zeta_m=e^{2\pi{\bm i}/m}$ throughout this section. We begin with a useful result due to Newman \cite[Theorem 1]{Newman}.
	
	\begin{lemma}\label{Lem. Newman}
		Let $m\ge 2$ be a positive integer. Suppose that there exist non-zero integers $c_1,c_2,\cdots,c_l$ and integers $0\le e_1<e_2<\cdots<e_l\le m-1$ such that 
		$$c_1\zeta_m^{e_1}+c_2\zeta_m^{e_2}+\cdots+c_l\zeta_m^{e_l}=0.$$
		Then $l\ge l(n)$, where $l(m)$ denotes the smallest prime factor of $m$. 
	\end{lemma}
	
	Recall that $N=p^m$ with $n=\varphi(n)=p^{m-1}(p-1)$. Now we prove our theorem.
	
	{\noindent\bfseries Proof of Theorem \ref{Thm. A}}. Throughout this proof, we fix a generator $g$ of $(\mathbb{Z}/N\mathbb{Z})^{\times}$. Given any generator $\chi$ of the dual group $\widehat{(\mathbb{Z}/N\mathbb{Z})^{\times}}$, we set $\chi(g)=\zeta_n^s$ for some $s\in\mathbb{Z}$ with $\gcd(s,n)=1$.
	
	Let $k$ be a positive divisor of $n$ with $n=kd$. For any integers $i,j\in[0,d-1]$, we first consider the $(i+1,j+1)$-entry of $A_k(\chi)$. One can verify that 
	\begin{align}\label{Eq. a in the proof of Thm. A}
	  A_k(\chi)(i+1,j+1)
&=G_N(\chi^{ki+kj})\notag\\
&=\sum_{x\in(\mathbb{Z}/N\mathbb{Z})^{\times}}\chi^{k(i+j)}(x)\zeta_N^x\notag\\
&=\sum_{0\le r\le n-1}\chi^{k(i+j)}(g^r)\zeta_N^{g^r}\notag\\
&=\sum_{0\le r\le n-1}\zeta_n^{sk(i+j)}\zeta_N^{g^r}\notag\\
&=\sum_{0\le r\le n-1}\zeta_d^{s(i+j)r}\zeta_N^{g^r}.
	\end{align}
	Noting that 
	$$[0,n-1]\cap\mathbb{Z}=\left\{a+bd: 0\le a\le d-1, 0\le b\le k-1\right\}$$
	and recalling that 
	$$H_d=\left\{x^d: x\in(\mathbb{Z}/N\mathbb{Z})^{\times}\right\}=\left\{g^{db}: 0\le b\le k-1\right\}\leq\left(\mathbb{Z}/N\mathbb{Z}\right)^{\times},$$
	by (\ref{Eq. a in the proof of Thm. A}) we obtain 
	\begin{align}\label{Eq. b in the proof of Thm. A}
		 A_k(\chi)(i+1,j+1)
&=\sum_{a=0}^{d-1}\sum_{b=0}^{k-1}\zeta_d^{s(i+j)(a+bd)}\zeta_N^{g^{a+bd}}\notag\\
&=\sum_{a=0}^{d-1}\zeta_d^{sa(i+j)}\sum_{b=0}^{k-1}\zeta_N^{g^{a+bd}}\notag\\
&=\sum_{a=0}^{d-1}\zeta_d^{sa(i+j)}\eta_N^{(g^a)}(k), 
	\end{align} 
	where 
	\begin{equation}\label{Eq. definition of eta N(k)}
		\eta_N^{(g^a)}(k)=\sum_{b=0}^{k-1}\zeta_N^{g^{a+bd}}=\sum_{x\in g^aH_d}\zeta_N^x
	\end{equation}
	is a Gauss period. Applying (\ref{Eq. b in the proof of Thm. A}), we obtain the matrix decomposition 
	\begin{equation}\label{Eq. c in the proof of Thm. A}
		A_k(\chi)=VDV,
	\end{equation}
	where $V$ is a $d\times d$ symmetric matrix defined by 
	$$V(i+1,j+1)=\zeta_d^{sij}$$
	for any integers $i,j\in[0,d-1]$, and 
	$$D=\diag\left(\eta_N^{(g^0)}(k),\eta_N^{(g^1)}(k),\cdots,\eta_N^{(g^{d-1})}(k)\right)$$
	is a diagonal matrix.  In addition, for any integers $i,j\in[0,d-1]$, it is easy to verify that 
	$$V^2(i+1,j+1)=\sum_{r=0}^{d-1}\zeta_d^{sr(i+j)}=\begin{cases}
		d  &  \mbox{if}\ i+j\equiv 0\pmod d,\\
		0  &  \mbox{otherwise}.
	\end{cases}$$
	From this we obtain 
	\begin{equation}\label{Eq. det of the square of V}
		(\det V)^2=(-1)^{\lfloor\frac{d-1}{2}\rfloor}\cdot d^d.
	\end{equation}
	
	Next we divide our remaining proof into two cases.
	
	{\bfseries Case 1}. $k\equiv 0\pmod p$. 
	
	In this case, clearly $m\ge2$. Let 
	$$C_p=\left\{1+p^{m-1}t\mod{N\mathbb{Z}}: 0\le t\le p-1\right\}\leq (\mathbb{Z}/N\mathbb{Z})^{\times}$$ 
	be the unique subgroup of order $p$. Since $k\equiv 0\pmod p$, we have $C_p\leq H_d$ by the Sylow theorem. Let 
	$$H_d=\bigcup_{1\le r\le k/p}y_rC_p$$
	be the coset decomposition of $H_d$ respect to $C_p$. Since
	$$\sum_{t\in\mathbb{Z}/p\mathbb{Z}}\zeta_p^{yt}=0$$
	for any $y\in (\mathbb{Z}/N\mathbb{Z})^{\times}$, one can verify that 
	\begin{align*}
		\eta_N^{(g^0)}(k)=\sum_{x\in H_d}\zeta_N^x
&=\sum_{r=1}^{k/p}\sum_{z\in C_p}\zeta_N^{y_rz}\\
&=\sum_{r=1}^{k/p}\sum_{t\in \mathbb{Z}/p\mathbb{Z}}\zeta_N^{y_r(1+p^{m-1}t)}\\
&=\sum_{r=1}^{k/p}\zeta_N^{y_r}\sum_{t\in\mathbb{Z}/p\mathbb{Z}}\zeta_p^{y_rt}.\\
&=0.
	\end{align*}
Now applying this to (\ref{Eq. c in the proof of Thm. A}), we immediately obtain 
$$\det A_k(\chi)=\left(\det V\right)^2\cdot \prod_{a=0}^{d-1}\eta_N^{(g^a)}(k)=0.$$
This completes the proof of (i).

{\bfseries Case 2}. $k\not\equiv 0\pmod p$. 

It is known that the Galois group 
$$\Gal\left(\mathbb{Q}(\zeta_N)/\mathbb{Q}\right)=\left\{\sigma_b: b\in\left(\mathbb{Z}/N\mathbb{Z}\right)^{\times}\right\}\cong \left(\mathbb{Z}/N\mathbb{Z}\right)^{\times},$$
where the $\mathbb{Q}$-automorphism $\sigma_b$ is determined by $\sigma_b(\zeta_N)=\zeta_N^b$ for any $b\in(\mathbb{Z}/N\mathbb{Z})^{\times}$. 

Now we evaluate the Galois group $\Gal(\mathbb{Q}(\eta_N^{(g^0)}(k))/\mathbb{Q})$. Since $p\nmid k$ and $$n=\varphi(p^m)=p^{m-1}(p-1)\equiv 0\pmod k,$$
we clearly have $p-1\equiv 0\pmod k$. Suppose first that $1\le k<p-1$. Note that $2k<p$ and that 
$$bH_d\cap H_d=\emptyset$$
whenever $b\in (\mathbb{Z}/N\mathbb{Z})^{\times}\setminus H_d$. By Lemma \ref{Lem. Newman}, when $1\le k<p-1$ we have  
\begin{equation}\label{Eq. d in the proof of Thm. A}
	\sigma_b\left(\eta_N^{(g^0)}(k)\right)-\eta_N^{(g^0)}(k)=\sum_{x\in H_d}\zeta_N^{bx}-\sum_{x\in H_d}\zeta_N^x\neq 0
\end{equation}
for any $b\in (\mathbb{Z}/N\mathbb{Z})^{\times}\setminus H_d$. Suppose now $k=p-1$. Then $d=p^{m-1}$. Given any $b\in (\mathbb{Z}/N\mathbb{Z})^{\times}\setminus H_d$, choose integers $u_1,u_2,\cdots,u_{p-1}, v_1,v_2,\cdots,v_{p-1}\in[0,N-1]$ such that 
$$\{u_r \mod{N\mathbb{Z}}: 1\le r\le p-1\}=H_d\ \text{and}\ \{v_r \mod{N\mathbb{Z}}: 1\le r\le p-1\}=bH_d.$$
Let
$$g(T)=\sum_{1\le r\le p-1}T^{v_r}-\sum_{1\le r\le p-1}T^{u_r}\in\mathbb{Z}[T]$$
be a polynomial with $\deg(g)<N=p^m$. If $g(\zeta_N)=0$, then 
$$g(T)\equiv 0\pmod{\Phi_{p^m}(T)\mathbb{Z}[T]},$$
where 
$$\Phi_{p^m}(T)=1+T^d+T^{2d}+\cdots+T^{(p-1)d}$$
is the $p^m$-th cyclotomic polynomial. Hence, there exists a polynomial $h(T)\in\mathbb{Z}[T]$ such that 
$$g(T)=h(T)\Phi_{p^m}(T).$$
This clearly implies that $\deg(h)<d$. Thus, we can set 
$$h(T)=c_0+c_1T+\cdots+c_{d-1}T^{d-1}$$
for some integers $c_0,c_1,\cdots,c_{d-1}$. Applying this to the decomposition $g(T)=h(T)\Phi_{p^m}(T)$, we see that in the polynomial $g(T)$, the coefficients of the $p$ terms $T^r, T^{r+d},\cdots, T^{r+(p-1)d}$ are all equal to $c_r$ for any integer $r\in [0,d-1]$. This contradicts the structure of the polynomial $g(T)$. Hence $g(\zeta_N)\neq 0$. This actually shows that 
\begin{equation}\label{Eq. d1 in the proof of Thm. A}
	\sigma_b\left(\eta_N^{(g^0)}(p-1)\right)-\eta_N^{(g^0)}(p-1)=\sum_{x\in H_d}\zeta_N^{bx}-\sum_{x\in H_d}\zeta_N^x=g(\zeta_N)\neq 0
\end{equation}
for any $b\in (\mathbb{Z}/N\mathbb{Z})^{\times}\setminus H_d$. 

On the other hand, given any positive divisor $k$ of $n$ and any $b\in H_d$, we clearly have 
\begin{equation}\label{Eq. e in the proof of Thm. A}
	\sigma_b\left(\eta_N^{(g^0)}(k)\right)=\sum_{x\in H_d}\zeta_N^{bx}=\sum_{x\in H_d}\zeta_N^{x}=\eta_N^{(g^0)}(k).
\end{equation}
Assembling (\ref{Eq. d in the proof of Thm. A}), (\ref{Eq. d1 in the proof of Thm. A}) and (\ref{Eq. e in the proof of Thm. A}) gives 
\begin{equation}\label{Eq. f in the proof of Thm. A}
	\Gal\left(\mathbb{Q}(\zeta_N)/\mathbb{Q}(\eta_N^{(g^0)}(k))\right)=\left\{\sigma_b: b\in H_d\right\}\cong H_d,
\end{equation}
and hence 
\begin{equation}\label{Eq. g in the proof of Thm. A}
		\Gal\left(\mathbb{Q}(\eta_N^{(g^0)})/\mathbb{Q}\right)\cong (\mathbb{Z}/N\mathbb{Z})^{\times}/H_d.
\end{equation}
From this, we see that $\mathbb{Q}(\eta_N^{(g^0)}(k))$ is the unique intermediate field of the cyclic Galois extension  $\mathbb{Q}(\zeta_N)/\mathbb{Q}$ with $[\mathbb{Q}(\eta_N^{(g^0)}(k)): \mathbb{Q}]=d$. 

Finally, we compute $\det A_k(\chi)$. Combining (\ref{Eq. f in the proof of Thm. A}) with (\ref{Eq. g in the proof of Thm. A}) and observing that 
$$(\mathbb{Z}/N\mathbb{Z})^{\times}/H_d=\left\{g^aH_d: 0\le a\le d-1\right\},$$
 we see that the minimal polynomial of the algebraic integer $\eta_N^{(g^0)}(k)$ over $\mathbb{Q}$ is equal to
$$P_k(T)=\prod_{b\in (\mathbb{Z}/N\mathbb{Z})^{\times}/H_d}\left(T-\sigma_b(\eta_N^{(g^0)}(k))\right)=\prod_{a=0}^{d-1}\left(T-\eta_N^{(g^a)}(k)\right)\in\mathbb{Z}[T].$$
Thus, the constant term $y_N(k)$ of $P_k(T)$ is equal to
\begin{equation}\label{Eq. h in the proof of Thm. A}
	y_N(k)=(-1)^d\cdot \prod_{a=0}^{d-1}\eta_N^{(g^a)}(k)\in\mathbb{Z}.
\end{equation}
Combining (\ref{Eq. h in the proof of Thm. A}) and (\ref{Eq. det of the square of V}) with (\ref{Eq. c in the proof of Thm. A}), we obtain 
$$\det A_k(\chi)=\left(\det V\right)^2\cdot \prod_{a=0}^{d-1}\eta_N^{(g^a)}(k)=(-1)^{d+\lfloor\frac{d-1}{2}\rfloor}\cdot d^d\cdot y_N(k).$$

In view of the above, we have completed the proof. \qed

	\Ack\ This research was supported by the Natural Science Foundation of China (Grant Nos. 12671009 and 12201291) and the Natural Science Foundation of the Higher Education Institutions of Jiangsu Province (Grant No. 25KJB110010).

\end{document}